 \documentclass[sn-mathphys,Numbered,pdflatex,referee]{sn-jnl}

\usepackage{mathtools}
\usepackage{graphicx}%
\usepackage{multirow}%
\usepackage{amsmath,amssymb,amsfonts}%
\usepackage{amsthm}%
\usepackage{amssymb}
\usepackage{mathrsfs}%
\usepackage[title]{appendix}%
\usepackage{xcolor}%
\usepackage{textcomp}%
\usepackage{manyfoot}%
\usepackage{booktabs}%
\usepackage{algorithm}%
\usepackage{algorithmicx}%
\usepackage{algpseudocode}%
\usepackage{listings}%
\usepackage{tikz}
\usetikzlibrary{quantikz2}

\begin{document}

\title[Quantum Diophantine Solver]{A Quantum Diophantine Equation Solution Finder}


\author[1]{\fnm{Lara} \sur{Tatli}}
\email{lara.tatli@durham.ac.uk}

\author*[2]{\fnm{Paul} \sur{Stevenson}}\email{p.stevenson@surrey.ac.uk}


\affil[1]{\orgdiv{Department of Physics}, \orgname{University of Durham}, \orgaddress{\city{Durham}, \postcode{DH1 3LE}, \country{UK}}}

\affil[2]{\orgdiv{School of Mathematics and Physics}, \orgname{University of Surrey}, \orgaddress{\city{Guildford}, \state{Surrey}, \postcode{GU2 7XH},  \country{UK}}}


\abstract{Diophantine equations are multivariate equations, usually polynomial, in which only integer solutions are admitted.  A brute force method for finding solutions would be to systematically substitute possible integer solutions and check for equality.  

Grover's algorithm is a quantum search algorithm which can find marked indices in a list very efficiently.  By treating the indices as the integer variables in the diophantine equation, Grover's algorithm can be used to find solutions in brute force way more efficiently than classical methods.  We present an example for the simplest possible diophantine equation.}




\maketitle

\section{Introduction}\label{intro}

A Diophantine equation is an equation, typically polynomial with integer coefficients, in more than one integer variable.  A famous example occurs as Fermat's Last Theorem, which states that
\begin{equation}
    x^n+y^n=z^n
\end{equation}
has no solutions for $n\ge3$ where $n$, $x$, $y$, and $z$ are all natural numbers.  The simplest Diophantine equation is linear in two variables and is of the form
\begin{equation}\label{lineareq}
    ax+by=n,
\end{equation}
where $a$, $b$, and $n$ are given constants.  While this equation has well-known solutions, in many other cases, solutions are not known (see e.g. the regularly-updated paper by Grechuk keeping track of some open and solved problems \cite{grechuk2022diophantine}).  Seeking solutions to Diophantine solutions through numerical search is an established method, where searches can prove the existence of solutions where it is posited that none exist \cite{74138}.  

Here, we bring quantum computing to bear upon the search for Diophantine equation solutions, using Grover's algorithm \cite{10.1145/237814.237866} to look for solutions for the simplest linear equation of the form (\ref{lineareq}), with $a=b=1$ and $n=5$ arbitrarily chosen for definiteness, and as a simple example that can be encoded in a workable number of qubits on an available simulator.  While we are not aware of works explicitly solving Diophantine equations with a quantum search algorithm, we note recent work using Grover's algorithm to perform a series of basic arithmetic procedures \cite{roy_applying_2024}.  In our work we use a standard classical adder for arithmetic and use Grover for the search for equality.

\section{Grover's Algorithm as equation solution searcher} 
We give here a brief discussion of the principles of a quantum search algorithm, following the treatment in Nielsen and Chuang's textbook  \cite{nielsen_chuang_2010}.  The search algorithm generally searches through a search space of $N$ elements.  It is supposed that one can work at the level of the index of the elements such that if presented with the index, it is easy to check if it is the element sought.  This is the case in our example where checking if given numbers $x$ and $y$ are solutions of the given equation is straightforward by direct substitution and evaluation.

The algorithm uses an oracle, $\mathcal{O}$, which acts as 
\begin{equation}
\mathcal{O}|x\rangle|q\rangle\rightarrow|x\rangle|q\oplus f(x)\rangle.\label{eq:oracle}
\end{equation}
Here, $|x\rangle$ is a register of index qubits, and $|q\rangle$ is the oracle qubit.  $\oplus$ is addition modulo 2 and $f(x)$ is a function which returns $0$ if index $x$ is not a solution to the search problem, and $1$ if index $x$ is a solution.

If the oracle qubit is prepared in the state $|-\rangle=(|0\rangle-|1\rangle)/\sqrt{2}$ then the action of the oracle is

\begin{equation}
    \mathcal{O}|x\rangle\left(\frac{|0\rangle-|1\rangle}{\sqrt{2}}\right)\rightarrow(-1)^{f(x)}|x\rangle\left(\frac{|0\rangle-|1\rangle}{\sqrt{2}}\right),
\end{equation}
thus the action of the oracle marks out, with a phase change, components of the register state $|x\rangle$ which are solutions to the problem - i.e. have $f(x)=1$.  The full Grover algorithm then amplifies the states which have been marked, and suppresses the unmarked states, using a ``diffuser'' circuit.  The oracle-diffuser combination together constitute a single Grover iteration, and a total of $O(\sqrt{N/M})$ iterations are needed in general to have the solutions selected in the register with high probability, where $M$ is the number of solutions in the $N$-element space.  Note that the standard diffuser requires that valid solutions do not account for the majority of the solution space, but this is the usual condition for an interesting Diophantine equation.

The indexing register works in our case by having $2m$ qubits in which each half encodes one of the numbers $x$ and $y$.  The encoding is made directly in standard binary and we do not consider negative numbers.  Clearly the size of $m$ will determine the available integers in the search space, and one must apply ever more qubits to increase the size of the search space, though one benefits from an exponential increase in search space as the number of qubits increases linearly.

For this exploratory study, to find solutions to the equation $x+y=5$ we use a $2m=6$ qubit register $|x\rangle$ to encode two 3-bit numbers to add together.  The oracle performs the addition and checks the result against the desired solution. The details of the quantum adder we use is given in the next section. 

\section{Quantum Adder circuit}

A quantum adder capable of calculating the sum of two 3-qubit binary numbers was produced using Qiskit. The adder was designed in such a way that the registers storing the input numbers were not overwritten during the calculation, as is the case with e.g. ripple-carry adders \cite{orts_review_2020}.  Retaining the input numbers is useful for use in further calculation, though not vital in our case.

\begin{figure}[tb]
\centering
 \begin{quantikz}[column sep = small]
        \lstick{$x_0$} & & & & & & & & & & \ctrl{9} &&\ctrl{8}&&\ctrl{8}&&\\
        \lstick{$x_1$} & & & &  \ctrl{9} & &\ctrl{6} & & \ctrl{6} & & && && &&\\
        \lstick{$x_2$} & \ctrl{9}& & \ctrl{4}& & & & & & & && && &&\\
        \lstick{$y_0$} & & & & & & & & & & &\ctrl{6}&\ctrl{5}&& &\ctrl{5}\slice{C}&\\
        \lstick{$y_1$} & & & & & \ctrl{6} &\ctrl{3} & & & \ctrl{3}\slice{B}& && && &&\\
        \lstick{$y_2$} & &\ctrl{6}&\ctrl{1} & & & & & & & && && &&\\
        \lstick{$a_0$} & & & \targ{}\slice{A} & & & &\ctrl{4} &\ctrl{1} & \ctrl{1}& && && &&\\
        \lstick{$a_1$} & & & & & &\targ{}& &\targ{} & \targ{}& && &\ctrl{2}&\ctrl{1} &\ctrl{1}&\\
        \lstick{$s_0$} & & & & & & & & & & &&\targ{}&&\targ{} &\targ{}&\\
        \lstick{$s_1$} & & & & & & & & & & \targ{} &\targ{}& &\targ{}& &&\\
        \lstick{$s_2$}& & & & \targ{}& \targ{} & &\targ{} & & & && && &&\\
        \lstick{$s_3$}& \targ{}&\targ{}& & & & & & & & && && &&\\
        \lstick{$q_{12}$}& & & & & & & & & & & & & & & &
    \end{quantikz}
\caption{A diagram of the quantum adder with barriers included to visually indicate each section.}\label{adder}
\end{figure}
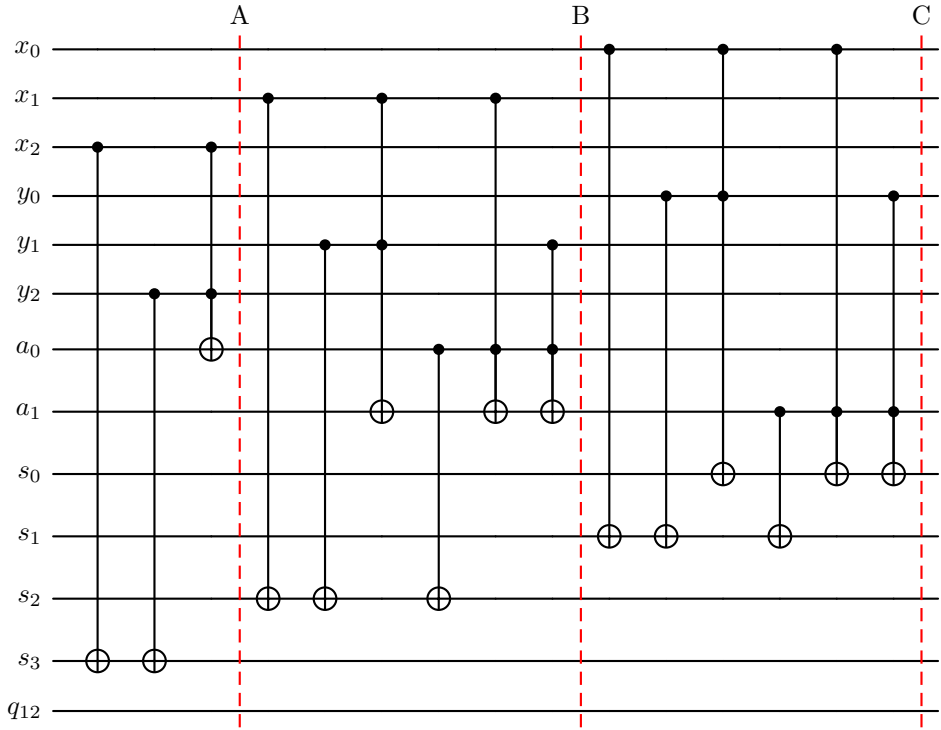

In this setup, shown in Fig. \ref{adder}, the first 3 qubits, $x_0$, $x_1$ and $x_2$, denote the binary digits representing a natural number $x$ in the format $x_0x_1x_2$, where $x_2$ is the least significant bit. In the same manner, qubits $y_0$, $y_1$ and $y_2$ denote the natural number $y$ in the format $y_0y_1y_2$. Qubits $a_0$ and $a_1$ represent ancillary qubits used to hold carry bits in the addition. Qubits $s_0$, $s_1$, $s_2$ and $s_3$ denote the solution to $x + y$ in the form $s_0s_1s_2s_3$, where $s_3$ is the least significant bit.   The figure shows all qubits in that are needed for the full Grover algorithm.  Qubit $q_{12}$ is the oracle qubit $|q\rangle$ as in equation (\ref{eq:oracle}).

The dividers labelled A, B, and C in the circuit help label different functional parts.  

In the section terminated by divider A, an addition operation is performed on the qubits representing the least significant bits $x_2$ and $y_2$ using two CNOT gates and one Tofolli gate, with the result stored in the qubit $s_3$ and the first carry bit stored in $a_0$.

In the section between dividers A and B, the qubits representing $x_1$, $y_1$, and the carry bit $a_0$ are added using three CNOT gates; the target is set to the sum digit $s_2$. Three Tofolli gates are used to compute the second carry bit, stored in $a_1$.

In the section between B and C, the sum digit $s_1$ is calculated using three CNOT gates acting on the qubits representing $x_0$, $y_0$, and the second carry bit in $a_1$. The final sum digit, $s_0$, is calculated using three Toffoli gates and takes into consideration the second carry bit.

In total, this adder employs 8 CNOT gates and 7 Toffoli gates collectively acting over 12 qubits.  In terms of scaling to larger registers, adding two $m$-bit numbers requires $4m$ qubits ($2m$ representing the numbers to be added, $m-1$ ancillary carry bits, and $m+1$ to represent the sum).  The number of gates is $3m-1$ CNOT gates and $3m-2$ Tofolli gates.

\section{Application of Grover's Algorithm}

In order to apply Grover's algorithm to solve a linear Diophantine equation $ax + by = n$ in the case $a=b=1$ and $n=5$, it is first necessary to apply a Hadamard gate to each of the qubits $|x_0\ldots x_2,y_0\ldots y_2\rangle$ encoding $x$ and $y$. This produces the initial superposition state with all possible solution strings present with equal amplitude.

\begin{figure}[tbh]
\centering
 \begin{quantikz}[column sep = small]
        \lstick{$x_0$} & & & & \\
        \lstick{$x_1$} & & & & \\
        \lstick{$x_2$} & & & & \\ 
        \lstick{$y_0$} & & & & \\
        \lstick{$y_1$} & & & & \\
        \lstick{$y_2$} & & & & \\
        \lstick{$a_0$} & & & & \\
        \lstick{$a_1$} & & & & \\
        \lstick{$s_0$} & \gate{X} & \ctrl{4} &\gate{X}&\\
        \lstick{$s_1$} & & \ctrl{3} & & \\
        \lstick{$s_{2}$} & \gate{X}&\ctrl{2} &\gate{X} & \\
        \lstick{$s_{3}$} & & \ctrl{1}& & \\
        \lstick{$q_{12}$} & & \targ{}& &
    \end{quantikz}
    \qquad\qquad
    \begin{quantikz}[column sep = small]
        \lstick{$x_0$} & \gate{H}&\gate{X} & \ctrl{12}& \gate{X} & \gate{H}& \\
        \lstick{$x_1$} & \gate{H}& \gate{X}& \ctrl{11}& \gate{X} & \gate{H}& \\
        \lstick{$x_2$} & \gate{H}& \gate{X}& \ctrl{10}& \gate{X} & \gate{H}& \\
        \lstick{$y_0$} & \gate{H}& \gate{X}& \ctrl{9}& \gate{X} & \gate{H}& \\
        \lstick{$y_1$} & \gate{H}& \gate{X}& \ctrl{8}& \gate{X} & \gate{H}& \\
        \lstick{$y_2$} & \gate{H}& \gate{X}& \ctrl{7}& \gate{X} & \gate{H}& \\
        \lstick{$a_0$} & & & &  & & \\
        \lstick{$a_1$} & & & &  & & \\
        \lstick{$s_0$}  & & & &  & & \\
        \lstick{$s_1$}  & & & &  & & \\
        \lstick{$s_{2}$} & & & &  & & \\
        \lstick{$s_{3}$} & & & &  & & \\
        \lstick{$q_{12}$}  & & &\targ{}&  & & \\
    \end{quantikz}
\caption{Left: Diagram of the query circuit and its inverse used for the oracle operation, $\mathcal{O}$, for the case $|s_0s_1s_2s_3\rangle$ = $|0101\rangle$. This circuit is run after the quantum adder circuit and is followed by the inverse quantum adder, forming a complete oracle.  Right: The diffuser circuit used to amplify the solution(s).}\label{oracle}
\end{figure}
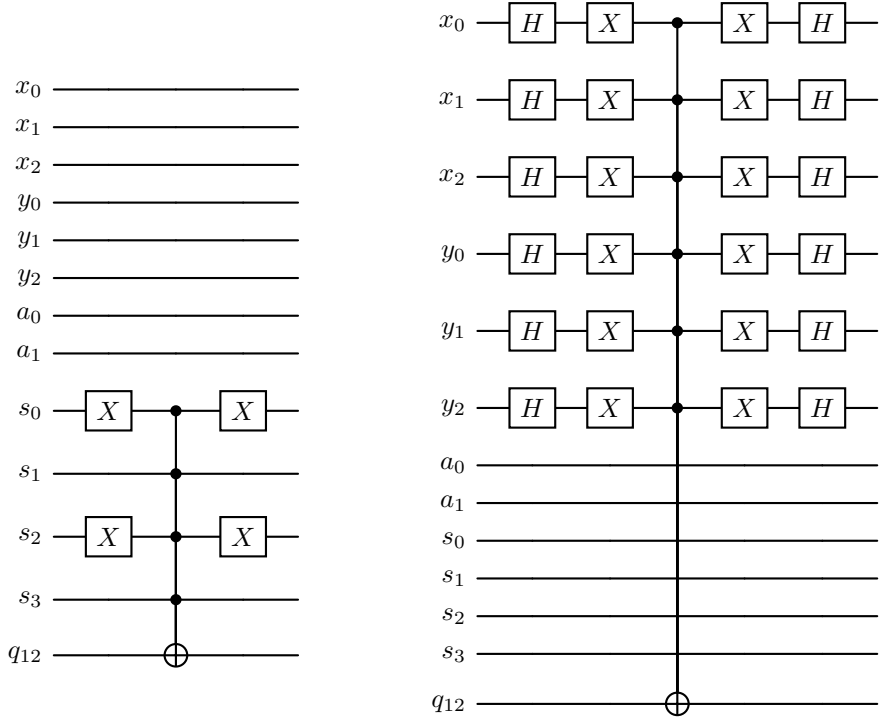

We then construct a quantum oracle capable of ``marking'' the solutions once queried.  This consists of the quantum adder and its inverse circuit with a query circuit in between which applies a phase shift of -1 to the solution qubits of the adder, if and only if, the solution is in the state $|s_0s_1s_2s_3\rangle$ = $|0101\rangle$. All other states are left unchanged. This is achieved using two X-gates and a multi-controlled Toffoli gate targeting $q_{12}$, configured to be in the $|-\rangle$ state prior to implementing Grover's algorithm. X-gates are re-applied to reverse the computation. The query circuit design used for this example is provided in the left-hand part of Fig. \ref{oracle}.

Each iteration of the oracle is followed by the circuit used for the diffusion operator, which by acting across the six qubits $|x_0\ldots x_2,y_0\ldots y_2\rangle$ amplifies states that sum to give the desired solution only. In this diffuser circuit, shown for our case in the right-hand part of Fig. \ref{oracle}, the combination of Hadamard and X-gates, in conjunction with a multi-controlled Toffoli gate, enable a phase change of -1 to be applied to the initial superposition state.  This completes one full iteration of the Grover algorithm.   After the desired number of algorithms, one would then perform a measurement on a real quantum computer, identically prepared through many repeated experiments, to build up a histogram of most probable outcomes corresponding to the sought solution(s).  In our present work, we simulate our circuit using a full quantum statevector, so present results in the next section by simply reading off the amplitudes of each register state.

\section{Implementation and Result}


\begin{figure}
    \begin{quantikz}
        \lstick{$x_0$} & \gate{H} & \gate[12]{Quantum Adder}& & \gate[12]{Quantum Adder\textsuperscript{\textdagger}} & \gate[13]{Diffuser} & \meter{}\\
        \lstick{$x_1$} & \gate{H} & & & & & \meter{} \\
        \lstick{$x_2$} & \gate{H} & & & & & \meter{} \\
        \lstick{$y_0$} & \gate{H} & & & & & \meter{}\\
        \lstick{$y_1$} & \gate{H} & & & & & \meter{}\\
        \lstick{$y_2$} & \gate{H} & & & & & \meter{}\\
        \lstick{$a_0$} & & & & & &\\
        \lstick{$a_1$} & & & & & &\\
        \lstick{$s_0$} & & & \gate[5]{Query}& & &\\
        \lstick{$s_1$} & & & & & &\\
        \lstick{$s_2$} & & & & & &\\
        \lstick{$s_3$} & & & & & &\\
        \lstick{$q_{12}$} & \gate{X} & \gate{H} & & & &\\
    \end{quantikz}
\caption{The complete circuit employing one Grover iteration.  The $\dagger$ symbol indicates Hermitian conjugate} \label{groverCircuit}
\end{figure}
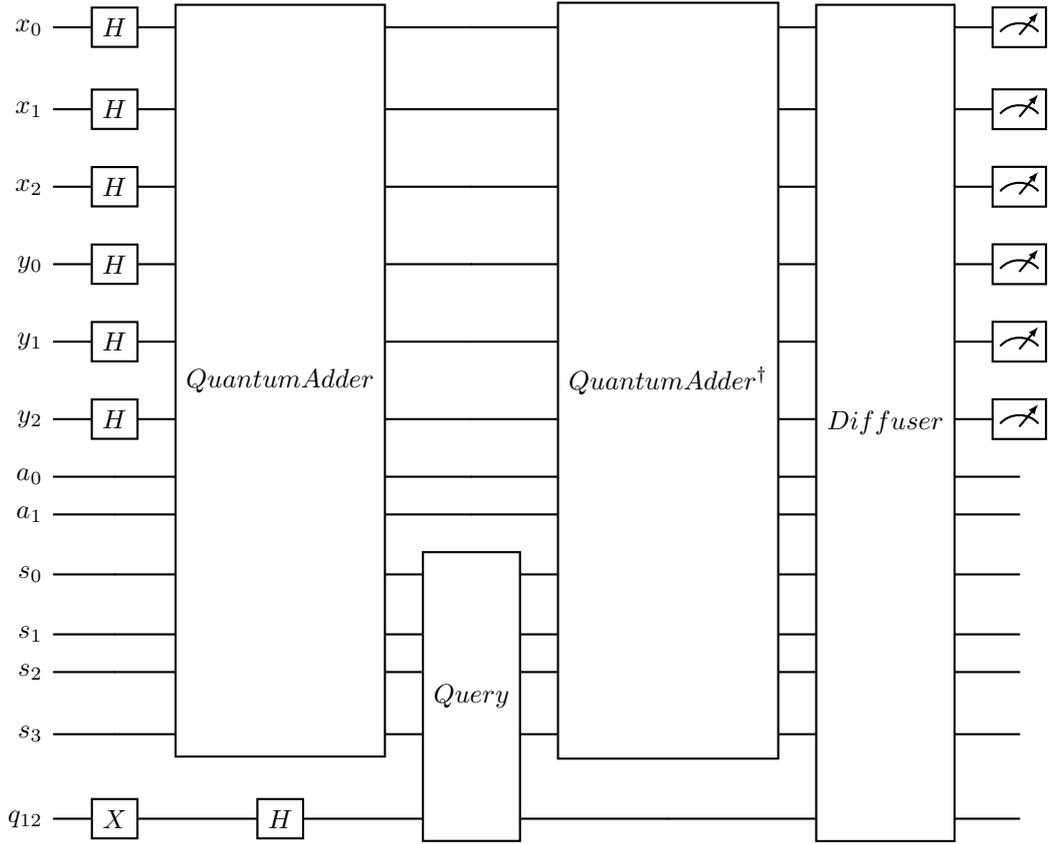

The full quantum circuit, including the Hadamards to initialize the superposition of the $x$ and $y$ register qubits and the $|-\rangle$ initialization of the oracle qubit, is shown in for one iteration in Fig \ref{groverCircuit}.  By running this full quantum circuit on BlueQubit's 34-qubit CPU statevector simulator, it is shown that two iterations of Grover's algorithm are sufficient to generate the full set of solutions to our simple Diophantine equation. 

\begin{figure}[tbh]
\centering
\includegraphics[width=1\textwidth]{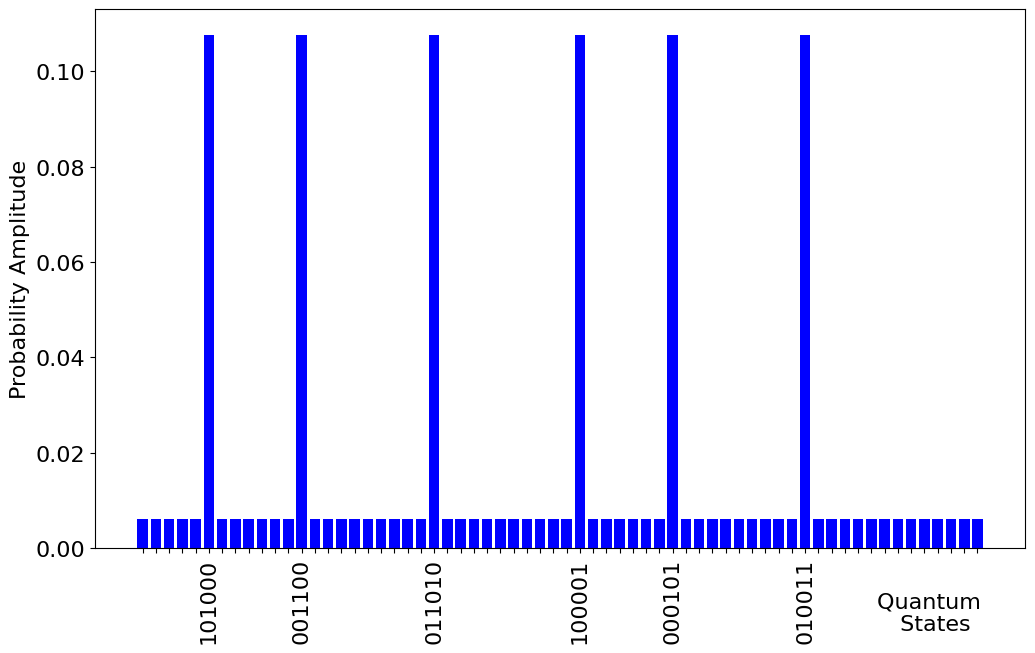}
\caption{Histogram for n=1 iterations.}\label{histogram_iterations1}
\end{figure}

\begin{figure}[tbh]
\centering
\includegraphics[width=1\textwidth]{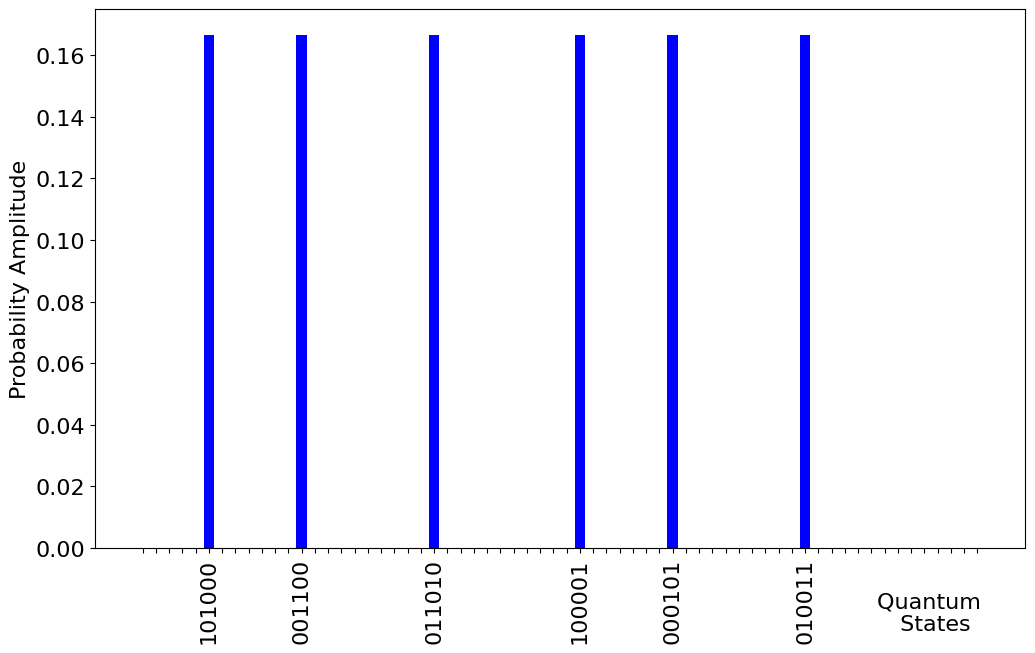}
\caption{Histogram for n=2 iterations - probabilities of incorrect solutions effectively become zero}\label{histogram_iterations2}
\end{figure}

The histogram displayed after one iteration is displayed in Fig. \ref{histogram_iterations1}; the histogram for two iterations is displayed in Fig. \ref{histogram_iterations2}. Note that the solution should be read from left to right, with the first three digits representing $x_0x_1x_2$ and the following digits $y_0y_1y_2$.

The solutions are seen to be correct solutions of the Diophantine equation $x+y=5$, and we tabulate them for clarity in Table \ref{tab:results}

\begin{table}[htb]
\begin{tabular}{cccccc}
\toprule
quantum state & $x$ (base 2) & $y$ (base 2) & $x$ (base 10) & $y$ (base 10) & $x+y$ (base 10) \\
\midrule
101000 & 101 & 000 & 5 & 0 & 5\\
001100 & 001 & 100 & 1 & 4 & 5\\
011010 & 011 & 010 & 3 & 2 & 5\\
100001 & 100 & 001 & 4 & 1 & 5\\
000101 & 000 & 101 & 0 & 5 & 5\\
010011 & 010 & 011 & 2 & 3 & 5\\
\bottomrule
\end{tabular}
\caption{Solution states picked out by Grover's algorithm in search for solutions to Diophantine equation $x+y=5$.}\label{tab:results}
\end{table}

We find that six iterations of Grover's algorithm are required to return to the probability distribution shown in Fig. \ref{histogram_iterations1}. 

\section{Conclusions}

Grover's algorithm can be implemented to search for solutions to simple linear Diophantine equations.  We have not attempted implementation on a real quantum computer, and the ability of our circuit to operate on noisy intermediate-scale quantum devices would need to be evaluated.  Nevertheless, further work could investigate more complicated Diophantine equations; for example, when $a$ and $b$ are not unity, or where the variables are raised to powers greater than unity.  In that case, more interesting unsolved cases, like those listed in Grechuk's paper \cite{grechuk2022diophantine} could be tackled.

Furthermore, we have not attempted to refine or optimize the algorithm, rather concentrating on a straightforward implementation.  Techniques to improve the Grover convergence \cite{abdulrahman_enhancing_2024} could be applied, while inclusion of a quantum counting approach \cite{aaronson_quantum_2019} would allow one to gain knowledge of how many Grover iterations should be applied in advance of performing each calculation.  For a more general Diophantine equation solver, such enhancements would be desirable.

\section*{Acknowledgements}
PDS acknowledges support from UK STFC under grant ST/W006472/1.
We acknowledge the use of IBM Quantum Lab in the early stages of this work.

\bibliography{pds-bib}

\end{document}